\documentclass[11pt,final,reqno]{article}
\usepackage{fullpage}
\usepackage{amsmath,amsfonts,amssymb,amsthm}
\usepackage{graphicx}
\usepackage{epstopdf}
\usepackage{color}
\usepackage{epsf}
\usepackage{verbatim}
\usepackage{enumitem}
\usepackage{hyperref}
\usepackage{subfig}
\usepackage[pagewise]{lineno}

\newtheorem{theorem}{Theorem}[section]
\newtheorem{lemma}[theorem]{Lemma}

\newtheorem{definition}[theorem]{Definition}
\newtheorem{conjecture}[theorem]{Conjecture}
\newtheorem{claim}[theorem]{Claim}
\newtheorem{fact}[theorem]{Fact}

\newcommand{\cT}{{\mathcal T}}
\newcommand{\cP}{{\mathcal P}}

\newcommand{\ep}{\varepsilon}

\newcommand{\floor}[1]{\left\lfloor#1\right\rfloor}
\newcommand{\ceiling}[1]{\left\lceil#1\right\rceil}

\def\beq{\begin{equation}}\def\eeq{\end{equation}}
\def\beqn{\begin{eqnarray}}\def\eeqn{\end{eqnarray}}

\def\eps{\varepsilon}
\def\({\mbox{$($}}\def\){\mbox{$)$}}

\let\OLDthebibliography\thebibliography
\renewcommand\thebibliography[1]{
  \OLDthebibliography{#1}
  \setlength{\parskip}{0pt}
  \setlength{\itemsep}{0pt plus 0.3ex}
}

\title{Ore-degree threshold for the square of a Hamiltonian cycle}
\author{Louis DeBiasio\thanks{Department of Mathematics, Miami University, Oxford, OH 45056 USA. E-mail address: debiasld@miamioh.edu.  Research supported in part by Simons Foundation Collaboration Grant \#402337.} \quad Safi Faizullah\thanks{Hewlett-Packard and Department of Computer Science, Rutgers University, Piscataway, NJ 08854, USA. E-mail address: safi@hp.com}  \quad Imdadullah Khan\thanks{Department of Computer Science, School of Science and Engineering, Lahore University of Management Sciences, Lahore 54792, Pakistan. E-mail address: imdad.khan@lums.edu.pk. Research supported in part by LUMS Faculty Startup Grant. Part of the work was done while the author was at Umm Al-Qura Univeristy, KSA.}}
\date{}

\begin{document}
\maketitle

\begin{abstract}
A classic theorem of Dirac from 1952 states that every graph with minimum degree at least $n/2$ contains a Hamiltonian cycle. In 1963, P\'osa conjectured that every graph with minimum degree at least $2n/3$ contains the square of a Hamiltonian cycle. 
In 1960, Ore relaxed the degree condition in the Dirac's theorem by proving that every graph with $\deg(u) + \deg(v) \geq n$ for every $uv \notin E(G)$ contains a Hamiltonian cycle. Recently, Ch\^au proved an Ore-type version of P\'osa's conjecture for graphs on $n\geq n_0$ vertices using the regularity--blow-up method; consequently the $n_0$ is very large (involving a tower function).  Here we present another proof that avoids the use of the regularity lemma.  Aside from the fact that our proof holds for much smaller $n_0$, we believe that our method of proof will be of independent interest.
\end{abstract}

\section{Introduction}

\subsection{Notation and Definitions}

Given a graph $G$, we denote the vertex set and edge set by $V(G)$ and $E(G)$ respectively, when the graph $G$ is clear by the context we refer to them as $V$ and $E$ respectively. When $uv\in E(G)$ we denote it by $u\sim v$ otherwise $u\not\sim v$. For a vertex $v\in V$, $N(v)$ is the set of neighbors of $v$ in $V$ and the degree of $v$ is $|N(v)|$ and we denote it by $\deg(v)=\deg_G(v)$. For $A\subseteq V(G)$, $N(v,A)$ is the set of neighbors of $v$ in $A$ and $\deg(v,A)$ is $|N(v,A)|$. 
 We denote by $\delta(G)$ the minimum degree over all vertices in $G$ and by  $\Delta(G)$ the maximum degree over all vertices in $G$.  We write $N(v_1,v_2,\dots, v_l)=\bigcap_{i=1}^{l}N(v_i)$ for the set of common neighbors of $v_1, v_2, \ldots, v_l$. Similarly, $N(v_1,v_2,\dots, v_l, A):=\bigcap_{i=1}^{l}N(v_i,A)$ and $deg(v_1,v_2,\ldots,v_l, A) := |N(v_1,v_2,\dots, v_l, A)|$.  We denote a cycle on $t$ vertices by $C_t$ and a path on $t$ vertices by $P_t$. When $G$ is a graph on $n$ vertices and $C_n\subseteq G$, we call $C_n$ a Hamiltonian cycle. A bipartite graph $G=(V,E)$, where $V = A \cup B$, $A\cap B=\emptyset$ will be denoted by $G(A,B)$.  The balanced complete $r$-partite graph with color classes of size $t$ is denoted by $K_r(t)$. For $A\subseteq V(G)$, $G[A]$ is the restriction of $G$ to $A$.
 When $A$ and $B$ are subsets of $V(G)$, we denote by $e(A,B)$ the number of edges of $G$ with one endpoint in $A$ and the other in $B$, and by $e(A)=|E(G[A])|$ the number of edges with both endpoints in $A$. Let $\delta(A,B)=\min_{v\in A} \deg(v, B)$.

For non-empty $A$ and $B$, $$d(A,B)=\frac{e(A,B)}{|A||B|}$$ is the \emph{density} of the graph between $A$ and $B$. We write $d(A)=2e(A)/|A|^2$. A graph $G$ on $n$ vertices is $\gamma$-\emph{dense} if it has at least $\gamma \binom{n}{2}$ edges. A bipartite graph $G(A,B)$ is $\gamma$-\emph{dense} if it contains at least $\gamma |A||B|$ edges.  Throughout the paper $\log$ denotes the base 2 logarithm.

\subsection{Powers of Cycles}

A classical result of Dirac \cite{DI} asserts that if $G$ is a graph on $n\geq 3$ vertices with $\delta(G)\geq n/2$, then $G$ contains a Hamiltonian cycle.  Note that when $n=2t$, Dirac's theorem implies that $G$ contains $t$ vertex disjoint copies of $K_2$.  In 1963, Corr\'adi and Hajnal \cite{CH} proved that if $G$ is a graph on $n=3t$ vertices with $\delta(G)\geq \frac{2n}{3}$, then $G$ contains $t$ vertex disjoint triangles.  Generalizing the Corr\'adi-Hajnal theorem, Erd\H{o}s conjectured \cite{E} and Hajnal and Szemer\'edi later proved \cite{HS} the following: 

\begin{theorem}[Hajnal-Szemer\'edi]\label{HSthm}
Let $G$ be a graph on $n=t(k+1)$ vertices. If $\delta(G)\geq \frac{kn}{k+1}$, then $G$ contains $t$ vertex disjoint copies of $K_{k+1}$.  
\end{theorem}

Finally in 1976, Bollobas and Eldridge \cite{BE}, and independently Catlin \cite{Cat}, made a conjecture which would generalize the Hajnal-Szemer\'edi theorem: If $G$ and $H$ are graphs on $n$ vertices with $\Delta(H)\leq k$ and $\delta(G)\geq \frac{kn-1}{k+1}$, then $H\subseteq G$.  While this conjecture is still open in general, we will only be interested in the $k=2$ case which was proved by Aigner and Brandt in 1993 \cite{AB}.
\begin{theorem}[Aigner-Brandt]\label{ABthm}
Let $G$ and $H$ be graphs on $n$ vertices.  If $\Delta(H)\leq 2$ and $\delta(G)\geq \frac{2n-1}{3}$, then $H\subseteq G$.
\end{theorem}
Note that all of these degree conditions are easily seen to be best possible.

Let $H$ be a graph with vertex set $V$.  The $k^{th}$ power of $H$, denoted $H^k$, is defined as follows: $V(H^k) = V$ and $uv\in E(H^k)$ if and only if the distance between $u$ and $v$ in $H$ is at most $k$. When $k=2$ we call $H^2$ the \emph{square of $H$}.  
Notice that $C_n^{k-1}$ contains $\floor{\frac{n}{k}}$ vertex disjoint copies of $K_k$.  Furthermore, notice that $C_n^2$ contains every graph $H$ on $n$ vertices with $\Delta(H)\leq 2$ (actually $P_n^2$ also has this property).  In 1963, P\'osa made a conjecture (see \cite{E}) that would significantly strengthen the Corr\'adi-Hajnal theorem (and retroactively Theorem \ref{ABthm}, see \cite{FK3}).
\begin{conjecture}[P\'osa]\label{posa}
Let $G$ be a graph on $n$ vertices. If $\delta(G)\geq\frac{2n}{3}$,
then $C_n^2\subseteq G$.
\end{conjecture}

After Erd\H{o}s' conjecture became the Hajnal-Szemer\'edi theorem, Seymour made a conjecture in 1974 \cite{SE} which generalizes P\'osa's conjecture to handle all values of $k$ (note that for $k\geq 4$, this does not generalize the Bollob\'as-Eldridge, Catlin conjecture).
\begin{conjecture}[Seymour]\label{sey}
Let $G$ be a graph on $n$ vertices. If $\delta(G)\geq\frac{kn}{k+1}$,
then $C_n^{k}\subseteq G$.
\end{conjecture}

Starting in the 90's a substantial amount of progress was made on these conjectures. Jacobson (unpublished) first established that the square of a Hamiltonian cycle can be found in any graph $G$ given that $\delta(G)\geq 5n/6$. Later Faudree, Gould, Jacobson and Schelp \cite{F2} improved the result, showing that the square of a Hamiltonian cycle can be found if $\delta(G)\geq(3/4+\eps)n$. The same authors further relaxed the degree condition to $\delta(G)\geq3n/4$. Fan and H\"aggkvist lowered the bound first in \cite{FH} to $\delta(G) \geq 5n/7$ and then in \cite{FK1} to $\delta(G)\geq(17n+9)/24$. Faudree, Gould and Jacobson \cite{F1} further lowered the minimum degree condition to $\delta(G)\geq7n/10$. Then Fan and Kierstead \cite{FK2} achieved the almost optimal $\delta(G)\geq\left(\frac23+\eps\right)n$. They also proved in \cite{FK3} that $\delta(G)\geq(2n-1)/3$ is sufficient for the existence of the square of a Hamiltonian {\em path}. Finally, they proved in \cite{FK4} that if $\delta(G)\geq2n/3$ and $G$ contains the square of a cycle with length greater than $2n/3$, then $G$ contains square of a Hamiltonian cycle. 

Regarding Conjecture \ref{sey}, Faudree {\it et al} \cite{F2} proved that for any $\eps > 0$ and positive integer $k$ there is a $C$ such that if $G$ is a graph on $n\geq C$ vertices with
$\delta(G)\geq\left(\frac{2k-1}{2k}+\eps\right)n,$
then $G$ contains the $k^{th}$ power of a Hamiltonian cycle.

Using the regularity--blow-up method first in \cite{KSS4} Koml\'os, S\'ark\"ozy and Szemer\'edi proved Conjecture \ref{sey} in asymptotic form,
then in \cite{KSS2} and \cite{KSS5} they proved both conjectures for $n\geq n_0$. The proofs used the regularity lemma \cite{SZ1}, the blow-up lemma \cite{KSS3, KSS6}, and the Hajnal-Szemer\'edi theorem \cite{HS}. Since the proofs used the regularity lemma the resulting $n_0$ is very large (it involves a tower function).  A new proof of P\'{o}sa's conjecture was given by Levitt, S\'ark\"ozy and Szemer\'edi \cite{LSS} which avoided the use of the regularity lemma and thus significantly decreased the value of $n_0$.  An explicit bound on $n_0$ was determined by Ch\^au, DeBiasio, and Kierstead in \cite{CDK};  however, for small $n_0$ the conjecture is still open. Finally, Jamshed and Szemer\'edi \cite{JSz} gave a new proof of  Seymour's conjecture that avoided the use of the regularity lemma. 

\subsection{Ore-type generalizations of Dirac-type results}

For a pair of non-adjacent vertices $(u,v)$, the value of $\deg(u)+\deg(v)$ is called the Ore-degree of $(u,v)$. We denote by $\delta_2(G)$ the minimum Ore-degree over all non-adjacent pairs of vertices in $G$.  In 1960, Ore \cite{O} proved that if $G$ is graph on $n\geq 3$ vertices with $\delta_2(G) \geq n$, then $G$ contains a Hamiltonian cycle.  Since any graph with $\delta(G)\geq \frac{n}{2}$ satisfies $\delta_2(G)\geq n$, Ore's theorem strengthens Dirac's theorem.  Inspired by this, researchers have sought to generalize minimum degree (``Dirac-type") conditions to Ore-type degree conditions; for a survey of such results see \cite{KKY}.  

Two important examples of Ore-type results are the following generalizations of Theorem \ref{HSthm} and \ref{ABthm}.  

\begin{theorem}[Kierstead-Kostochka \cite{KK}]\label{KKthm}
Let $G$ be a graph on $n=t(k+1)$ vertices.  If $\delta_2(G)\geq \frac{2kn}{k+1}-1$, then $G$ contains $t$ vertex disjoint copies of $K_{k+1}$.
\end{theorem}

\begin{theorem}[Kostochka-Yu \cite{KY}]\label{KYthm}
Let $G$ and $H$ be graphs on $n$ vertices.  If $\Delta(H)\leq 2$ and $\delta_2(G)\geq \frac{4n}{3}-1$, then $H\subseteq G$.
\end{theorem}

A natural Ore-type generalization of P\'osa's conjecture suggests that if $\delta_2(G) \geq \frac{4n}{3}$, then $C_n^2\subseteq G$.  It turns out that this natural generalization is not quite true as Ch\^au \cite{C} gave a construction of a graph $G$ for which $\delta_2(G)=\frac{4n}{3}$, but $G$ does not contain the square of a Hamiltonian cycle. However, in the same paper, Ch\^au uses the regularity--blow-up method to prove that if $G$ is a graph on $n\geq n_0$ vertices with $\delta_2(G)>\frac{4n}{3}$, then $C_n^2\subseteq G$.  In fact, he is able to give an even more refined degree condition:

\begin{theorem}[Ch\^au]\label{Cthm}  Let $G$ be a graph on $n$ vertices. If $\delta_2(G)\geq \frac{4n-1}{3}$ and
\begin{enumerate}
\item $\delta(G)\leq \frac{n}{3}+2$, then $P_n^2\subseteq G$. 

\item $\delta(G)>\frac{n}{3}+2$, then there exists $n_0$ such that if $n\geq n_0$, then $C_n^2\subseteq G$.
\end{enumerate}
\end{theorem}
\noindent
(See \cite{C}, Proposition 9.1 for an explanation of why this result actually implies Theorem \ref{KYthm} and the $k=2$ case of Theorem \ref{KKthm} for sufficiently large $n$ despite the fact that $\frac{4n-1}{3}>\frac{4n}{3}-1$.)

\begin{figure}[ht]
\centering
\subfloat[$\delta_2(G)=\frac{4n}{3}$ and $\delta(G)=\frac{n}{3}+2$]{
\makebox[.4\textwidth]{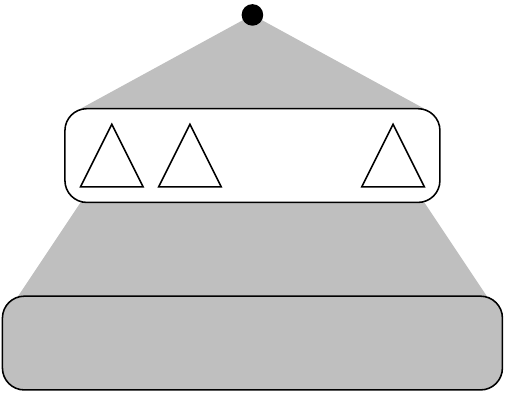}}
\hspace{1in}
\subfloat[$\delta_2(G)=\frac{4n-2}{3}$]{
\makebox[.4\textwidth]{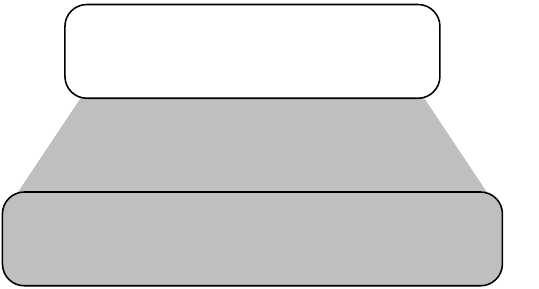}}
\caption[]{Examples showing the tightness of the degree conditions in Theorem \ref{main} }
\label{Examples}
\end{figure}

One of the purposes of this paper is to present another proof of Theorem \ref{Cthm}.(ii) which avoids the use of the regularity lemma, thus resulting in a much smaller value of $n_0$.

\begin{theorem}\label{main}  There exists $n_0$ such that if $G$ is a graph on $n\geq n_0$ vertices with 
\begin{equation} \label{minDeg}
\delta_2(G) \geq \frac{4n-1}{3} \mbox{ and } \delta(G) > \frac{n}{3} + 2,
\end{equation}
then $C_n^2\subseteq G$.
\end{theorem}

Aside from lowering the bound on $n_0$, we believe that the techniques used in this paper are of independent interest and can have more applications. In particular,  our proof provides a simpler template for approaching the following Ore-type version of Conjecture \ref{sey}.

\begin{conjecture}[Ch\^au]\label{oreseymour}
Let $G$ be a graph on $n$ vertices.  If $\delta_2(G)\geq \frac{2kn-1}{k+1}$ and $\delta(G)>\frac{(k-1)n}{k+1}+2$, then $C_n^k\subseteq G$.
\end{conjecture}

\subsection{Outline of the Proof}

As is common in these types of problems, our proof is divided into extremal and non-extremal cases.  The extremal conditions will resemble the properties found in Figure \ref{Examples}; either there is a vertex close to smallest possible degree, or there is a set of size approximately $n/3$ with very few edges.  We formally define the extremal conditions below.

\begin{definition}[Extremal]
Let $0<\alpha\ll\frac{1}{3}$ and let $G$ be a graph on $n$ vertices.
\begin{enumerate}
\item We say that $G$ satisfies extremal condition 1 with parameter $\alpha$ if there exists $v\in V(G)$ such that $\deg(v) < (\frac{1}{3}+\alpha)n$.

\item We say that $G$ satisfies extremal condition 2 with parameter $\alpha$ if there exist disjoint sets $A_1, A_2$ such that for $i=1,2$, $|A_i|\geq (1/3-\alpha)n$ and $d(A_i)<\alpha$.

\item We say that $G$ satisfies extremal condition 3 with parameter $\alpha$ if there exists a set $A_1$ such that $|A_1|\geq (1/3-\alpha)n$, $d(A_1)<\alpha$, and for all $A_2\subseteq V(G)\setminus A_1$ with $|A_2|\geq (1/3-\alpha)n$, $d(A_2)\geq \alpha$.
\end{enumerate}
\end{definition}

%
%
%

\begin{definition}[Non-extremal]
Let $0<\alpha\ll\frac{1}{3}$.  If $G$ does not satisfy extremal condition 1,2, and 3 with parameter $\alpha$, then we say $G$ is not $\alpha$-extremal.  Specifically, this implies that $\delta(G)\geq (1/3+\alpha)n$ and for all $A\subseteq V(G)$ with $|A|\geq (1/3-\alpha)n$, $d(A)\geq \alpha$.
\end{definition}

These extremal cases are dealt with in \cite{C} without the use of the regularity lemma; however, the blow-up lemma is used in multiple cases.  
In Section \ref{EXsection} we provide an alternate argument which can be used in \cite{C} to replace each use of the blow-up lemma.  


The non-extremal case is where our proof differs most significantly from \cite{C} and is the main focus of our paper.  We avoid the use of the regularity lemma, the blow-up lemma, and Theorem \ref{KKthm} by instead using Erd\H{o}s-Stone type results to cover all but a small fraction of the vertex set with disjoint balanced complete tripartite graphs of size about $\log n$. Then we prove a new connecting lemma which allows us to connect the complete tripartite graphs by square paths.  Aside from any leftover vertices, we have a nearly spanning structure which contains a square cycle and is quite robust  in the sense that most of the vertices are in complete tripartite graphs of size about $\log n$.  Finally, we take advantage of the robustness of our structure by inserting the leftover vertices in such a way that the resulting structure contains the square of a Hamiltonian cycle.  All of this will be made precise in Section \ref{nonex_section}.

\section{Extremal case}\label{EXsection}

In \cite{C}, the extremal cases are handled with very detailed, yet elementary arguments -- with one exception.  In certain cases of \cite{C}, the problem is reduced to finding the square of a Hamiltonian cycle in a balanced tripartite graph where each pair is nearly complete, with the exception of a small number of vertices which still satisfy some minimum degree condition.  Here Ch\^au uses the fact that these very dense pairs are $(\ep, \delta)$-super regular so the blow-up lemma can be applied to show that the desired square cycle exists.  However, these dense pairs have a property which is far stronger than the property of being $(\ep, \delta)$-super regular. Thus, our goal in this section is simply to provide an elementary argument which could be used to replace all of the uses of the blow-up lemma in the extremal cases of \cite{C}. Note that we will not reproduce the proof found in \cite{C}, as we are only providing a minor diversion to the conclusion of certain cases of the argument.

\begin{lemma}\label{cleanExtCaseLemma1}
Let $0<\alpha' \ll 1$ and let $H$ be a balanced tripartite graph on $3m= n \geq n_0$ vertices with $V(H)$ partitioned as $A_1,A_2,A_3$.  If for all $i\neq j$, $\delta(A_i ,A_j) \geq (1-\alpha')m$, then we can cover $V(H)$ by disjoint triangles.
\end{lemma}
\begin{proof}
We first find a perfect matching $M_1$ between $A_1$ and $A_2$ by an application of the K\"onig-Hall theorem. Then we find a perfect matching between $M_1$ and $A_3$, such that $e=xy\in M_1$ is matched with a vertex $z \in N(x,y,A_3)$. For any edge $e=xy \in M_1$ we have $\deg(x,y,A_3) \geq (1-2\alpha')m$, therefore, by K\"onig-Hall theorem there exists a perfect matching between $M_1$ and $A_3$ as desired. 
\end{proof}

\begin{lemma}\label{cleanExtCaseLemma2}
Let $0<\alpha' \ll 1$  and let $H$ be a balanced tripartite graph on $3m= n \geq n_0$ vertices with $V(H)$ partitioned as $A_1,A_2,A_3$.  If $T = \{t_1,t_2,\ldots,t_m\}$ is a triangle cover of $V(H)$ and if for all $i\neq j$,  $\delta(A_i, A_j) \geq (1-\alpha')m$, then $H$ contains the square of a Hamiltonian cycle.  Furthermore, $H$ contains the square of a Hamiltonian path which starts with $t_1$ and ends with $t_m$.
\end{lemma}

\begin{proof}
Let $t=(x_1,x_2,x_3)$ and $t' = (y_1,y_2,y_3)$ be any two triangles in $T$ such that $x_i,y_i \in A_i$. We say that $t$ {\em precedes} $t'$, if $x_i$ is adjacent to $y_1,\ldots,y_{i-1}$ for $2 \leq i \leq 3$ (if $t$ {\em precedes} $t'$, then $x_1x_2x_3y_1y_2y_3$ is a square-path). We say that $\{t, t'\}$ is a {\em good} pair, if $t$ precedes $t'$ and $t'$ precedes $t$. By the degree conditions above, any $t_i\in T$ makes a {\em good} pair with at least $(1-\sqrt{\alpha'}) m$ other triangles in $T$.

Make an auxiliary graph $H'$ over $T$ such that each triangle $t_i\in T$ is adjacent to the triangle $t_j$ if and only if $\{t_i, t_j\}$ is a good pair. By the above observation we clearly have $\delta(H') > m/2$, hence by the Dirac's theorem there is a Hamiltonian cycle in $H'$.  Also since $\delta(H') > m/2$, $H'$ is Hamiltonian connected and thus there is a Hamiltonian path in $H'$ which starts with $t_1$ and ends with $t_m$.
It is easy to see that this Hamiltonian cycle (path) in $H'$ corresponds to the square of a Hamiltonian cycle (path) in $H$. 
\end{proof}

Finally we arrive at the main lemma which can be used to replace the use of the blow-up lemma in the extremal cases of \cite{C}.

\begin{lemma}
Let $0<\alpha' \ll \beta\ll \gamma\ll 1$ and let $H$ be a balanced tripartite graph on $3m= n \geq n_0$ vertices with $V(H)$ partitioned as $A_1,A_2,A_3$.  If for all $i\neq j$, there are at least $(1-\beta)m$ vertices in $A_i$ with at least $(1-\alpha')m$ neighbors in $A_j$ and $\delta(A_i, A_j)\geq \gamma m$, then $H$ contains the square of a Hamiltonian cycle.  Furthermore, if we specify two edges $u_1u_2$ and $u_{3m-1}u_{3m}$ such that for all $u\in \{u_1, u_2, u_{3m-1}, u_{3m}\}$, $\deg(u, A_j)\geq (1-\alpha')m$, then $H$ contains the square of a Hamiltonian path $P=u_1u_2\dots u_{3m-1}u_{3m}$.
\end{lemma}

\begin{proof}
Call a vertex $u$ in $A_i$ {\em bad} if $u$ has less than $(1-\alpha')m$ neighbors in $A_j$ for some $j\neq i$.  By the hypothesis, there are at most $2\beta m$ bad vertices in each $A_i$.
Now with a simple greedy procedure, for each bad vertex $u\in A_1$ we find a triangle $t_2 = (b_1,b_2,b_3)$, such that $b_1=u$ and $b_2$ and $b_3$ are typical (not bad) vertices in $A_2$ and $A_3$. We find two more similar triangles $t_1 = (a_1,a_2,a_3)$ and $t_3 = (c_1,c_2,c_3)$, such that $a_2\in N(u)$, $a_3\in N(u,b_2)$ and $c_1\in N(b_2,b_3)$, $c_2\in N(b_3)$. Clearly $a_1a_2a_3b_1b_2b_3c_1c_2c_3$ is a square path. We replace these three triangles with an exceptional triangle $(d_1,d_2,d_3)$ with one vertex each in $A_1$, $A_2$ and $A_3$, such that for $1\leq i\leq 3$, $d_i$ is connected to common neighbors of $a_i$ and $c_i$. By the fact that $a_i$ and $c_i$ are not bad vertices every $d_i$ has at least $(1-3\alpha')m$ neighbors in both of the other two sets.  We similarly make an exceptional triangle for each of the remaining bad vertices.  Since the total number of bad vertices is at most $6\beta m$ and the minimum degree is $\gamma m\gg 6\beta m$, this greedy procedure can be easily carried out.  In the remaining parts of $A_1$, $A_2$, and $A_3$ by Lemma \ref{cleanExtCaseLemma1} we find a triangle cover and add all the exceptional triangles to the cover.  Then by Lemma \ref{cleanExtCaseLemma2}, we find the square of a Hamiltonian cycle. 

Now suppose $u_1u_2$ and $u_{3m-1}u_{3m}$ are given edges such that for all $u\in \{u_1, u_2, u_{3m-1}, u_{3m}\}$, $\deg(u, A_j)\geq (1-\alpha')m$ for all $A_j$ such that $u\notin A_j$.  We make $t_1 = (u_1,u_2,u_3)$ and $t_m = (u_{3m-2},u_{3m-1},u_{3m})$ such that $u_3$ is a typical vertex in $N(u_1,u_2)$ and $u_{3m-2}$ is a typical vertex in $N(u_{3m-1},u_{3m})$.  Now by applying Lemma  \ref{cleanExtCaseLemma1} we find a triangle cover and add all the exceptional triangles to the cover. Then by Lemma \ref{cleanExtCaseLemma2}, we find the square of a Hamiltonian path which starts with $t_1$ and ends with $t_m$. 
\end{proof}

\section{Non-extremal case}\label{nonex_section}

Before we give an overview of the non-extremal case, it would be helpful to have some idea of how the non-extremal case is proved in \cite{C} (which is a generalization of the arguments in  \cite{KSS2}, \cite{KSS4}, \cite{KSS5}). Suppose $G$ is a non-extremal graph on $n$ vertices ($n$ sufficiently large) with $\delta_2(G)\geq \frac{4n}{3}$.  Using the regularity lemma and Theorem \ref{KKthm}, one can show that $G$ contains a set of disjoint balanced $4$-partite and $3$-partite graphs spanning almost all of $G$ each having size $\Omega(n)$.  Each of these multi-partite graphs $H$ has the property that every pair of color classes forms a suitably dense psuedorandom bipartite graph, so by applying the blow-up lemma, one obtains an almost spanning square path in $H$.  If we connect these multi-partite graphs together with square paths before applying the blow-up lemma, one will obtain an almost spanning square path of $G$.  Finally the remaining vertices need to somehow be inserted, which is an elementary, but detailed argument.  

We are able to avoid the regularity--blow-up method by showing that for sufficiently large $n$ (but nowhere near as big as needed for the regularity lemma), $G$ can be partitioned into disjoint balanced complete tripartite graphs spanning almost all of $G$, each having size $\Omega(\log n)$; we call this ``the cover" and it is built in Section \ref{sec:cover}.  Since the tripartite graphs are complete, we do not have to apply the blow-up lemma; if we go around a complete tripartite graph picking vertices from each of the color classes sequentially we get a square-path.  Next we must prove a Connecting Lemma which allows us to connect the tripartite graphs by short square-paths giving us a ``cycle of cliques"; this is done in Section \ref{sec:connecting}.  At the end of this process there will be a few leftover vertices which need to be inserted; this is done in Section \ref{sec:inserting}.



Here is the statement of the non-extremal case (notice that in the non-extremal case we are able to slightly relax the Ore-degree condition).

\begin{theorem}[Non-extremal case]\label{thm:nonextremal}
For all $0<\ep\ll \alpha\ll 1$ there exists $n_0$ such that if $G$ is a graph on $n\geq n_0$ vertices with 
\beq \label{minDegree_nonext} 
\delta_2(G)\geq (\frac{4}{3}-2\ep)n 
\eeq 
and $G$ is not $\alpha$-extremal, then $C_n^2\subseteq G$.
\end{theorem}

Given a graph $G$ with $\delta_2(G)\geq (\frac{4}{3}-2\ep)n$, let $L$ be the set of vertices in $G$ with degree less than $(2/3-\ep)n$.  We say that the vertices in $L$ are \emph{low-degree vertices} and the vertices in $V(G)\setminus L$ are \emph{high-degree vertices}.  Note that by the definition of $\delta_2(G)$, $G[L]$ is a clique.

\subsection{The Cover}\label{sec:cover}

In order to cover most of the vertices in $G$ with complete tripartite graphs as mentioned above, we will need quantitative versions of some classical results in extremal graph theory.
\subsubsection{Lemmas}\label{sec:lemmas}

%
%
%

\begin{fact}\label{denseSubgraph}
Let $0<d,\gamma<1$.  If $G(A,B)$ is a $(d+2\gamma)$-dense bipartite graph, then there must be at least $\gamma|B|$ vertices in $B$ for which the degree in $A$ is at least $(d+\gamma)|A|$.
\end{fact}

\begin{proof}
Indeed, otherwise the total number of edges would be less than
$$(d+\gamma) |A|\cdot|B| + \gamma |B|\cdot|A| = (d+2\gamma) |A||B|$$ a contradiction to the fact that $G(A,B)$ is $(d+2\gamma)$-dense. 
\end{proof}

\begin{lemma}\label{rpartVolArg}
Let $0<c,\gamma<1/3$, $s=\floor{c\log n}$, and let $G$ be a graph on $n\geq n_0$ vertices with $K:=K_3(s)=(A_1,A_2,A_3)\subseteq G$.  If $B\subseteq V(G)\setminus V(K)$ with $|B|\geq \gamma n$ and $d(B,K)\geq \frac{2}{3} + 2\gamma$, then $G[B\cup V(K)]$ contains a $K':=K_4(\gamma s) = (A_1',A_2',A_3',B')$, where $A_i'\subset A_i$ and $B'\subset B$.
\end{lemma}

\begin{proof}
Let $B_1=\{v\in B: \deg(v, K)\geq (\frac{2}{3} + \gamma)|V(K)|\}$. By Fact \ref{denseSubgraph} we have $|B_1|\geq \gamma |B|\geq \gamma^2 n$.  By the degree condition, each vertex in $B_1$ has at least $\gamma s$ neighbors in each $A_i$.  There are at most $2^{|V(K)|}=2^{3s} \le n^{3c}$ different possible neighborhoods, so by averaging there must be a neighborhood that appears for a set $B_2$ of at least $\frac{|B_1|}{n^{3c}} \geq \frac{\gamma^2n }{n^{3c}}=\gamma^2 n^{(1-3c)}$ vertices of $B_1$. Selecting an appropriate subset $B'$ of $B_2$, we get the desired complete $K_4(\gamma s)$. 
\end{proof}

We need a version of the Erd\H{o}s-Stone theorem where we have control of the parameters.  While there are a sequence of improvements by Bollob\'as-Erd\H{o}s, Bollob\'as-Erd\H{o}s-Simonovits, and Bollob\'as-Kohayakawa (to name a few), we will state a version due to Nikiforov \cite{Nik1} which gives an explicit lower bound on $n$.  

\begin{lemma}[\cite{Nik1}]\label{BES}
Let $c$ and $n$ be such that $0<c<1$ and $n\geq \exp(64/c^3)$, and let $G$ be a graph on $n$ vertices.  If $e(G)\geq (\frac{1}{2}+c)\frac{n^2}{2}$, then $G$ contains $K_{3}(s)$ where $s=\floor{\frac{c^3}{64}\log n}$.
\end{lemma}

Finally, we need a simple fact which allows us to translate our Ore-degree condition into an appropriate edge density condition.

\begin{fact}\label{OreEdges}
Let $0<d<1$ and let $G$ be a graph on $n\geq 2$ vertices. If $\delta_2(G) \geq 2dn$, then $e(G) \geq d\frac{n^2}{2}$.
\end{fact}

\begin{proof}
Define $\gamma$ so that $e(G)=\gamma\binom{n}{2}$ and suppose $\delta_2(G) \geq 2d(n-1)$.  We have
\begin{align*}
2d(n-1)(1-\gamma)\binom{n}{2}\leq \sum_{\{u,v\}\not\in E(G)}(\deg(u)+\deg(v))&=\sum_{v\in V(G)}\deg(v)(n-1-\deg(v))\\
&=2\gamma\binom{n}{2}(n-1)-\sum_{v\in V(G)}(\deg(v))^2\\
&\leq 2\gamma\binom{n}{2}(n-1)-\frac{1}{n}\left(\sum_{v\in V(G)}\deg(v)\right)^2\\
&= 2\gamma\binom{n}{2}(n-1)-\frac{4\gamma^2\binom{n}{2}^2}{n}.
\end{align*}
Dividing both sides by $2(n-1)(1-\gamma)\binom{n}{2}$ gives $\gamma\geq d$, and thus $\delta_2(G) \geq 2d(n-1)$ implies $e(G)\geq d\binom{n}{2}$.  Thus if $\delta_2(G)\geq 2dn=2\frac{dn}{n-1}(n-1)$, then $e(G)\geq \frac{dn}{n-1}\binom{n}{2}=d\frac{n^2}{2}$ as stated.
\end{proof}

\subsubsection{Building the cover}

\begin{definition}[Tripartite Cover]
Let $s, n'\in \mathbb{R}^+$. A $(s, n')$ tripartite cover is a collection $\cT$ of vertex disjoint copies of $K_3(t_i)=:T_i$ with $\floor{s}\leq t_i\leq \ceiling{2s}$ such that $|\bigcup_{T_i\in \cT}V(T_i)|\geq n'$.
\end{definition}

Note that in the following lemma we do not assume that $G$ is non-extremal.

\begin{lemma}[Cover Lemma]\label{lem:cover}
For all $0<\ep\ll \eta\ll 1$,
there exists $n_0$ and $c>0$ such that if $G$ is a graph on $n\geq n_0$ vertices with $\delta_2(G)\geq (\frac{4}{3}-2\ep)n$, then 
$G$ contains a $(c\log n, (1-\eta)n)$ tripartite cover.
\end{lemma}

\begin{proof}
Set $c_0=\frac{\eta^6}{64}$ and $t_0=\floor{c_0\log n}$.
By (\ref{minDegree_nonext}) and Fact \ref{OreEdges} we have $e(G) \geq (\frac{2}{3}-\ep)\frac{n^2}{2}$. 
We repeatedly apply Lemma \ref{BES} (with $c=\eta^2$) to find complete tripartite graphs with each color class of size $t_0$ until the remaining graph contains no copy of $K_3(t_0)$.  Let $\cT$ be the collection of tripartite graphs obtained in this way, and let $U=V(G)\setminus V(\cT)$, where $V(\cT)=\bigcup_{T\in \cT}V(T)$.  If $|U| < \eta n$, then we are done, so suppose $|U|\geq \eta n$.

Set $U_0=U$, $\cT_0=\cT$, and for $i\geq 0$ set $c_i=\eta^{2i}c_0=\eta^{2i+6}/64$ and $t_i=\floor{c_i\log n}$. 

\begin{claim}\label{enlargecover}
If $|U_i|\geq \eta n$, then either
\begin{enumerate}
\item $G[U_i]$ contains $K:=K_3(t_i)$, in which case we reset $U_i:=U_i\setminus V(K)$ and $\cT_{i}:=\cT_i \cup K$ or

\item $G[U_i]$ does not contain a copy of $K_3(t_i)$, in which case there exists a cover $\cT_{i+1}$ such that $|V(\cT_{i+1})|\geq |V(\cT_i)|+\eta^4n$ and every color class in the cover has size between $t_{i+1}$ and $2t_{i+1}$.
\end{enumerate}
\end{claim}

If Claim \ref{enlargecover} holds, then for some $j\leq \frac{1}{\eta^4}$, we must have $|U_j|<\eta n$ (as at least $\eta^4 n$ vertices are added to the cover before we increase the index).  Note that $c_j\geq \eta^{2j}\eta^6/64\geq \eta^{\frac{2}{\eta^4}+6}/64=:c$ and thus the balanced complete tripartite graphs in $\cT_j$ have parts of size between $\floor{c\log n}$ and $\ceiling{2c\log n}$ as required.  We now finish the proof of the cover lemma by proving Claim \ref{enlargecover}.

\begin{proof}

Let $0\leq i\leq \frac{1}{\eta^4}$ and suppose $G[U_i]$ does not contain a copy of $K_3(t_i)$.  In this case by Lemma \ref{BES} 
\beq \label{Udensity}
d(U_i)<(\frac{1}{2}+\eta^{(2i+6)/3})\leq (\frac{1}{2}+\eta^2).   
\eeq

Start by setting $Z=\emptyset$.  We will consider each $T\in \cT_i$ one by one.  If $d(U_i,T) < (\frac{2}{3}+6\eta^2)$, then consider the next element of $\cT_i$.  If $d(U_i,T) \geq (\frac{2}{3}+6\eta^2)$, then by Lemma \ref{rpartVolArg} there exists $K_{4}(3\eta^2 t_i)$ in $G[U_i\cup V(T)]$, which can be split into four copies of $K_3(\eta^2 t_i)=K_3(t_{i+1})$.  Move the used vertices from $U_i$ into $Z$ and reset $U_i:=U_i\setminus Z$.  Let $\cT_i'$ be the set of $3$-partite graphs in $\cT_i$ for which the procedure succeeded.  If $|\cT_i'|\geq \eta^2\frac{n}{3t_i}$, then we will have increased the cover by at least $3\eta^2 t_i\cdot\eta^2\frac{n}{3t_i} =\eta^4 n$.  If $|U_i|<\eta n$ or we have increased the cover by $\eta^4n$, we partition each color class into parts of size at least $t_{i+1}$ (which implies that all parts have size at most $2t_{i+1}$).

So suppose we have increased the cover by less than $\eta^4n$ and we still have $|U_i|\geq \eta n$.  In this case we have $|\cT_i'|<\eta^2\frac{n}{3t_i}$ which implies 
\beq \label{T'}
|V(\cT_i')\cup Z|= |\cT_i'|(3t_i+3\eta^2t_i)<3(1+\eta^2)t_i\cdot\eta^2\frac{n}{3t_i}<2\eta^2n.  
\eeq

For every $T\in \cT_i\setminus \cT_i'$, we have 
\begin{align} 
e(U_i, T) \leq
(\frac{2}{3}+6\eta^2)|V(T)||U_i|   \label{UtoT}.
\end{align}

Now by \eqref{T'} and \eqref{UtoT} we have 
\begin{align}
e(U_i, V(\cT_i)\cup Z)=e(U_i,V(\cT_i')\cup Z) + e(U_i, V(\cT_i\setminus \cT_i'))
&\leq 2\eta^2n |U_i| + (\frac{2}{3}+6\eta^2)|V(\cT_i)||U_i|\notag \\
&\leq (\frac{2}{3}|V(\cT_i)| + 8\eta^2 n)|U_i|\label{UtoOutside}
\end{align}

Recall that $G[L]$ (the graph induced by the low-degree vertices) induces a clique and since $G[U_i]$ contains no $K_3(t_i)$, we have $|L\cap U_i|< 3t_i<\ep |U_i|<\ep n$.  Also note that $n=|U_i|+|Z|+|V(\cT_i)|\geq |U_i|+|V(\cT_i)|$.  Now we get  
\begin{align}
e(U_i, V(\cT_i)\cup Z)=\sum_{v\in U_i} \deg(v) - 2 e(U_i) &\stackrel{\eqref{Udensity}}{\geq} \sum_{v\in U_i\setminus L} (\frac{2}{3}-\eps)n - 2(\frac{1}{2}+\eta^2)\frac{|U_i|^2}{2} \notag\\
&\geq ((1-\ep)(\frac{2}{3}-\eps)n - (\frac{1}{2}+\eta^2)|U_i|)|U_i|\notag\\
&\geq (\frac{2}{3}|V(\cT_i)| + \frac{2}{3}|U_i|- 2\eps n - (\frac{1}{2}+\eta^2)|U_i|)|U_i|\notag\\
&\geq (\frac{2}{3}|V(\cT_i)| + \frac{1}{6}|U_i|- 2\eps n -\eta^2|U_i|)|U_i|\notag\\
&\geq (\frac{2}{3}|V(\cT_i)| + \frac{1}{6}\eta n- 2\eps n -\eta^2|U_i|)|U_i|
\label{UtoOutside_lower}
\end{align}

By \eqref{UtoOutside} and \eqref{UtoOutside_lower}, we have $(\frac{2}{3}|V(\cT_i)|+\frac{1}{7}\eta n)|U_i|\leq e(U_i, V(\cT_i)\cup Z)\leq (\frac{2}{3}|V(\cT_i)| + 8\eta^2 n)|U_i|$, a contradiction.
\end{proof}

\end{proof}

\subsection{Connecting}\label{sec:connecting}

In this section we will prove that if $G$ is non-extremal, then we can find a short square path between any two disjoint ordered edges provided that each edge consists of high degree vertices or has a triangle in the common neighborhood of the endpoints.  This Lemma will then be used to connect the tripartite graphs coming from Lemma \ref{lem:cover}.

\subsubsection{Connecting ordered edges}

First note the following simple fact.

\begin{fact}\label{4or5edges}
Given disjoint triangles $T$ and $T'$, either there exists an ordering of the vertices of $T=x_3x_1x_2$ and $T'=y_1y_2y_3$ such that $x_3x_1x_2y_1y_2y_3$ is a square path, or there exist vertices $x_1,x_2\in T$ and $y_1,y_2\in T'$ such that $x_1\not\sim y_1$ and $x_2\not\sim y_2$. 
\end{fact}

\begin{proof}
One very special case of a result of Faudree and Schelp \cite{FS} says that in every $2$-coloring of $K_{3,3}$ there is either a red path on $4$ vertices or a blue path on $4$ vertices (this special case is easily checked).  Applying this to the induced bipartite graph between $T$ and $T'$ implies that (with the appropriate labeling of the vertices) $x_1y_1$, $y_1x_2$, and $x_2y_2$ are either all edges or all non-edges; the former implies that $x_3x_1x_2y_1y_2y_3$ is a square path and the latter implies that $x_1\not\sim y_1$ and $x_2\not\sim y_2$.
\end{proof}

\begin{lemma}[Connecting Lemma]\label{generalconnecting}
Let $\frac{1}{n_0}\ll \ep\ll \alpha\ll 1$ and let $G$ be a graph on $n\geq n_0$ vertices with $\delta_2(G)\geq (\frac{4}{3}-4\ep)n$ such that $G$ is not $\alpha$-extremal.  For all distinct $u_1,u_2,v_1,v_2\in V(G)$ with $u_1u_2, v_1v_2\in E(G)$, if 
\begin{enumerate}
\item $\deg(u_1), \deg(u_2) \geq (\frac{2}{3}-2\ep)n$ ~or~ there exists a triangle $T$ in $G[N(u_1,u_2)\setminus \{v_1,v_2\}]$, and 
\item $\deg(v_1), \deg(v_2) \geq (\frac{2}{3}-2\ep)n$ ~or~ there exists a triangle $T'$ in $G[N(v_1,v_2)\setminus \{u_1,u_2\}]$,
\end{enumerate}
then there exists $Q\subseteq G-\{u_1,u_2,v_1,v_2\}$ such that  $P=u_1u_2Qv_1v_2$ is a square path with $|V(Q)|\leq 18$.

\end{lemma}


%
%
%

\begin{proof}  Suppose first that there exists a triangle $T$ in $G[N(u_1,u_2)\setminus \{v_1,v_2\}]$ and there exists a triangle $T'$ in $G[N(v_1,v_2)\setminus \{u_1,u_2\}]$; let $G'=G-T- T'-\{u_1,u_2,v_1,v_2\}$.  If $T$ and $T'$ are vertex disjoint, then by Fact \ref{4or5edges}, we either immediately find a square path from $T$ to $T'$ or there exist two disjoint non-adjacent pairs of vertices in $T\times T'$.  Let $(x_i,y_i)$ and $(x_j,y_j)$ be two such pairs and define $C_{i,j}=\{v\in V(G'): \deg(v, \{x_i,y_i,x_j,y_j\})=4\}$.  Consider two disjoint non-edges $(x_i,y_i)$, $(x_j,y_j)$ such that $|C_{i,j}|$ is maximum.  We may label the vertices of $T$ as $x_1,x_2,x_3$ and the vertices of $T'$ as $y_1,y_2,y_3$ such that the disjoint non-edges which maximize $|C_{i,j}|$ are $(x_1,y_1)$ and $(x_2,y_2)$; i.e., $C_{i,j}=C_{1,2}$.  Let $A=\{v\in V(G): \deg(v,\{x_1,x_2\})=2$, $\deg(v, \{y_1, y_2\})=1\}$ and $B=\{v\in V(G): \deg(v,\{x_1,x_2\})=1$, $\deg(v, \{y_1, y_2\})=2\}$.  Set $C:=C_{1,2}$ and note that $A, B$ and $C$ are disjoint.  Since $x_1\not\sim y_1$ and $x_2\not\sim y_2$, we have 
\begin{equation}\label{x1x2y1y2}
\deg(\{x_1,y_1,x_2,y_2\})\geq 2\left(\frac{4}{3}-4\ep\right)n= \left(\frac{8}{3}-8\ep\right)n.
\end{equation}
Also we have
$$\deg(\{x_1,y_1,x_2,y_2\})\leq 4|C|+3(|A|+|B|)+2(n-|A|-|B|-|C|).$$
Together this gives
\begin{equation}\label{AB2C}
|A|+|B|+2|C|\geq \left(\frac{2}{3}-8\ep\right)n.
\end{equation}
If $T$ and $T'$ have an edge $e$ in common, then set $Q=e$ and note that $u_1u_2Qv_1v_2$ is the desired square path with $|Q|= 2$.  If $T$ and $T'$ have exactly one vertex in common, call it $z$.  Now, if there exist vertices $x\in V(T)$ and $y\in V(T')$, both distinct from $z$, such that $xy\in E(G)$, then set $Q=xzy$ and note that $u_1u_2Qv_1v_2$ is the desired square path with $|Q|= 3$.  So if $T$ and $T'$ have exactly one vertex in common, we can label the vertices of $T$ as $x_1,x_2,x_3$ and the vertices of $T'$ as $y_1,y_2,y_3$ such that $x_3=y_3$ and $x_1\not\sim y_1$ and $x_2\not\sim y_2$.  Letting $G'=G-\{u_1,u_2,v_1,v_2\}-T-T'$, we may define $A$, $B$, and $C$ as above, with \eqref{x1x2y1y2} and consequently \eqref{AB2C} holding.

Now suppose, without loss of generality, that there is no triangle in $G[N(u_1,u_2)\setminus \{v_1v_2\}]$ but there is a triangle $T'$ in $G[N(v_1,v_2)\setminus \{u_1,u_2\}]$.  We have $|N(u_1,u_2)\setminus L|\geq (\frac{1}{3}-8\ep)n-|L|>(\frac{1}{3}-\alpha)n$ and since we are not in the extremal case, we have an edge $x_1'x_2'\in G[N(u_1,u_2)\setminus L]-T'$.  If there is a triangle in $G[N(x_1',x_2')\setminus \{u_1,u_2,v_1,v_2\}]$, then we call it $T$ and proceed as in the first paragraph, noting that in this case $Q$ must begin with $x_1'x_2'$.  So suppose there is no triangle in $G[N(x_1',x_2')\setminus \{u_1,u_2,v_1,v_2\}]$. This implies $|N(x_1',x_2')\setminus L|\geq (\frac{1}{3}-8\ep)n-|L|>(\frac{1}{3}-\alpha)n$ and since we are not in the extremal case, we have an edge $x_1''x_2''\in G[N(x_1',x_2')\setminus L]$.  Again, if there is a triangle in $G[N(x_1'',x_2'')\setminus \{u_1,u_2,x_1',x_2',v_1,v_2\}]$, then we call it $T$ and proceed as in the first paragraph, noting that in this case $Q$ must begin with $x_1'x_2'x_1''x_2''$.  So suppose there is no triangle in $G[N(x_1'',x_2'')\setminus \{u_1,u_2,v_1,v_2,x_1',x_2'\}]$.  Suppose $x_1''$ or $x_2''$, say $x_1''$, has $3$ neighbors in $T\rq{}$. If there exists $i\in [2]$ such that $x_i'$ has at least one neighbor in $T'$, call it $y_1$ and let $y_2$ be a distinct vertex in $T'$, we may set $Q=x_{3-i}'x_i'x_1''y_1y_2$ and note that $u_1u_2Qv_1v_2$ is the desired square path with $|Q|=5$.  On the other hand, if $x_1''$ has $3$ neighbors in $T'$, but $x_1'$ and $x_2'$ have no neighbors in $T'$, then we set $x_1:=x_1'$ and $x_2:=x_2'$ and with $A$, $B$, $C$ defined as before, \eqref{x1x2y1y2} and consequently \eqref{AB2C} holds.  So suppose $x_1''$ and $x_2''$ each have less than $3$ neighbors in $T'$.  Either $x_1''$ and $x_2''$ have the same two neighbors in $T'$, say $y_1',y_2'$, in which case we set $Q=x_1'x_2'x_1''x_2''y_1'y_2'$, giving us the desired square path $u_1u_2Qv_1v_2$ with $|Q|=6$, or else there exists $y_1,y_2\in T'$ such that with $x_1:=x_1''$ and $x_2:=x_2''$, we have $x_1\not\sim y_1$ and $x_2\not\sim y_2$.  So with $A$, $B$, $C$ defined as before, \eqref{x1x2y1y2} and consequently \eqref{AB2C} holds.

Finally suppose that there is no triangle in $G[N(u_1,u_2)\setminus \{v_1v_2\}]$ and no triangle in $G[N(v_1,v_2)\setminus \{u_1,u_2\}]$; this implies that $|N(u_1,u_2)\cap L|< 3$ and $|N(v_1,v_2)\cap L|< 3$ .  So we have $|N(u_1,u_2)\setminus L|,|N(v_1,v_2)\setminus L|\geq (\frac{1}{3}-8\ep)n-|L|>(\frac{1}{3}-\alpha)n$ and since we are not in the extremal case, we have an edge $x_1'''x_2'''\in G[N(u_1,u_2)\setminus L]$ and a disjoint edge $y_1'''y_2'''\in G[N(v_1,v_2)\setminus L]$.  If either $x_1''',x_2'''$ has a triangle in their common neighborhood or $y_1''',y_2'''$ has a triangle in their common neighborhood, we proceed as in one of the previous cases, noting that we append $x_1'''x_2'''$ or $y_1'''y_2'''$ to the beginning or end of $Q$ respectively. So suppose not; in this case we set $x_i:=x_i'''$ and $y_i:=y_i'''$ for all $i\in [2]$.  Letting $G'=G-\{u_1,u_2,v_1,v_2,x_1,x_2,y_1,y_2\}$, we may define $A$, $B$, and $C$ as before and since $x_1,x_2,y_1,y_2\not\in L$, \eqref{x1x2y1y2} and consequently \eqref{AB2C} holds.

With the initial segments of the square path now in place, suppose there exists a square path $Q'$ having at most 10 vertices starting with either direction of $x_1x_2$ and ending with either direction of $y_1y_2$.  From the cases above, we may have to append at most 6 vertices to the beginning of $Q'$ and at most $2$ vertices to the end of $Q'$, giving us the desired square path $P=u_1u_2Qv_1v_2$ with $|Q|\leq 18$.

\begin{claim}\label{ABCempty}
If following conditions do not hold, then we have the desired square path $Q'$.
\begin{enumerate}
\item $d(C)=0$, $d(B,C)=0$, $d(A,C)=0$, and $d(A,B)=0$; and

\item for all $S\in \{A,B,C\}$, if $d(S)=0$, then $|S|\leq 1$.  In particular, $|C|\leq 1$.

\end{enumerate}
\end{claim}

\begin{proof}
\begin{enumerate}

\item From the definition of $A$, $B$, $C$, if any of the given edge sets were non-empty, we would immediately have a square path $Q'$. 
%
%
%

\item First suppose $|A|+|B|\leq \alpha n$.  Then by \eqref{AB2C}, we have
$$|C|\geq \frac{1}{2}\left(\frac{2}{3}-8\ep-\alpha\right)n\geq \left(\frac{1}{3}-\alpha\right)n.$$
Since we are not in the extremal case, we have $d(C)>0$ contradicting Claim \ref{ABCempty}(i).  So we may suppose 
\begin{equation}\label{A+B}
|A|+|B|>\alpha n
\end{equation}

Let $S\in \{A,B,C\}$ with $d(S)=0$ and suppose $|S|\geq 2$.  From Claim \ref{ABCempty}(i), we have $\deg(v, A\cup B\cup C)=0$ for all $v\in S$.  Let $v_1,v_2\in S$ and since $v_1\not\sim v_2$ we have $\deg(v_1)+\deg(v_2)\geq \left(\frac{4}{3}-4\ep\right)n$. But now we have the following contradiction
\begin{align*}
\left(\frac{4}{3}-4\ep\right)n\leq \deg(v_1)+\deg(v_2)&\leq 2(n-|A|-|B|-|C|)\\
&\leq 2n-(|A|+|B|+2|C|)-(|A|+|B|)\\
&\leq 2n-\left(\frac{2}{3}-8\ep\right)n-\alpha n ~~~\text{(by \eqref{AB2C} and \eqref{A+B})}\\
&<\left(\frac{4}{3}-4\ep\right)n.
\end{align*}

Note that since $d(C)=0$ by Claim \ref{ABCempty}(i), we have $|C|\leq 1$ by Claim \ref{ABCempty}(ii).

\end{enumerate}
\end{proof}

From this point we assume that the conditions of Claim \ref{ABCempty} hold, as otherwise we would be done.  By Claim \ref{ABCempty}(ii) and \eqref{AB2C}, we have
\begin{equation}\label{ABbig}
|A|+|B|\geq \left(\frac{2}{3}-8\ep\right)n-2|C|\geq \left(\frac{2}{3}-9\ep\right)n.
\end{equation}

We consider two cases based on the edge density of $A$ and $B$.

\noindent\textbf{Case 1:} $d(A)>0$ and $d(B)>0$

Let $a_1a_2\in G[A]$ and $b_1b_2\in G[B]$.  By Claim \ref{ABCempty}(i), $a_1\not\sim b_1$ and $a_2\not\sim b_2$.  Thus $|N(a_1, b_1)|,|N(a_2, b_2)|\geq \left(\frac{1}{3}-8\ep\right)n$.  Furthermore, by Claim \ref{ABCempty}(i), $|N(a, b)\cap (A\cup B\cup C)|=0$ for all $a\in A$, $b\in B$.  Combining this with \eqref{ABbig} gives  
\begin{align*}
|N(a_1,b_1,a_2,b_2)|\geq 2\left(\frac{1}{3}-8\ep\right)n-(n-|A\cup B\cup C|)&\geq \left(\frac{2}{3}-16\ep\right)n-\left(\frac{1}{3}+9\ep\right)n \\
&>\left(\frac{1}{3}-\alpha\right)n.
\end{align*}
Since we are not in the extremal case, we have an edge $d_1d_2\in G[N(a_1,b_1,a_2,b_2)]$ which gives the desired square path $Q'=x_1x_2a_1a_2d_1d_2b_1b_2y_1y_2$ with $|Q'|\leq 10$.

\noindent
\textbf{Case 2:} $d(A)=0$ or $d(B)=0$.

If $d(B)=0$, then by Claim \ref{ABCempty}(ii), $|B|\leq 1$ and since $|C|\leq 1$, \eqref{AB2C} implies 
\beq \label{Abig}
|A|\geq \left(\frac{2}{3}-8\ep\right)n-|B|-2|C|\geq \left(\frac{2}{3}-9\ep\right)n.
\eeq
Likewise if $d(A)=0$, then 
\beq \label{Bbig}
|B|\geq \left(\frac{2}{3}-8\ep\right)n-|A|-2|C|\geq \left(\frac{2}{3}-9\ep\right)n.
\eeq

Recall that so far we have constructed square paths $u_1u_2\dots \{x_1,x_2\}$ and $\{y_1,y_2\}\dots v_1v_2$ and we are attempting to construct $Q'$ to complete the square path from some ordering of $x_1x_2$ to some ordering of $y_1y_2$.  If $x_1x_2$ is in a triangle $T$ in the common neighborhood of the two vertices preceding $x_1x_2$ on $u_1u_2\dots x_1x_2$, we say that $x_1x_2$ is \emph{supported by $T$}.  Likewise, if $y_1y_2$ is in a triangle $T'$ in the common neighborhood of the two vertices following $y_1y_2$ on $y_1y_2\dots v_1v_2$, we say that $y_1y_2$ is \emph{supported by $T'$}.  This distinction is important, because if say $x_1x_2$ is not supported by $T$, then based on the initial construction this implies that $x_1,x_2$ are high degree vertices and have no triangle in their common neighborhood.

\textbf{Case 2.1:} $d(A)=0$ and $x_1x_2$ is supported by $T$, or $d(B)=0$ and $y_1y_2$ is supported by $T'$. 

Without loss of generality, suppose $d(B)=0$ and $y_1y_2$ is supported by $T'$. By \eqref{Abig} and $\delta(G)\geq (1/3+\alpha)n$, we have $$|A\cap N(y_3)|\geq (\alpha - 9\ep)n\geq 3.$$  
First suppose there exists $i\in \{1,2\}$ such that $y_3\sim x_i$.  Let $a\in N(y_3)\cap A$ and let $j\in \{1,2\}$ such that $a\sim y_j$ (by definition of $A$).  Thus $Q'=x_{3-i}x_iay_3y_jy_{3-j}$ is the desired square path with $|Q'|\leq 6$.  So suppose that $y_3\not\sim x_1$ and $y_3\not\sim x_2$ (if $x_1x_2$ is supported by $T$, this of course implies that we are in the case where $T$ and $T'$ do not have a vertex in common).  Since $|N(y_3)\cap A|\geq 3$, there exists $j\in \{1,2\}$ such that $y_j\sim a$ and $y_j\sim a'$ for distinct $a,a'\in N(y_3)\cap A$.  Furthermore, since $a,a'\in A$, we have $\{x_1,x_2,y_3,y_j\}\subseteq N(a)\cap N(a')$, but since $|C|\leq 1$ the non-adjacent pairs $(x_j, y_j)$ and $(x_{3-j}, y_3)$ along with $a,a'$ contradict the maximality of $|C|$.

\textbf{Case 2.2:}  $d(A)=0$ and $x_1x_2$ is not supported by $T$, or $d(B)=0$ and $y_1y_2$ is not supported by $T'$.

Without loss of generality, suppose $d(A)=0$ and $x_1x_2$ is not supported by $T$. By \eqref{Bbig}, we have $|B|\geq (\frac{2}{3}-9\ep)n$.  By the case we have $|N(x_1,x_2)|\geq (\frac{1}{3}-4\ep)n$ and there exists $z\in N(x_1,x_2)\setminus L$.  Since $\deg(z)\geq (\frac{2}{3}-2\ep)n$ and $|B|\geq (\frac{2}{3}-9\ep)n$, we have $|N(z)\cap B|\geq (\frac{1}{3}-11\ep)n>(\frac{1}{3}-\alpha)n$.  As $G$ is not $\alpha$-extremal, there exists $b_1b_2\in E(G[N(z)\cap B])$.  Since $b_1\in B$, there exists $i\in\{1,2\}$ such that $b_1\sim x_i$.  Thus $Q'=x_{3-i}x_izb_1b_2y_1y_2$ is the desired square path with $|Q'|\leq 7$.

\end{proof}

\subsubsection{Connecting complete tripartite graphs}

\begin{definition}[Connected tripartite cover]
Let $q, s ,n'\in \mathbb{R}^+$.  A $(q, s, n')$ connected tripartite cover is a $(s, n')$ tripartite cover $\{K^0, \dots, K^{m-1}\}$ together with a collection of $m$ vertex disjoint square paths $\{P_0, \dots, P_{m-1}\}$ such that for all $0\leq i\leq m-1$, $P_i=u_1u_2u_3Qv_1v_2v_3$ is a square path from $K^i=(U_1,U_2,U_3)$ to $K^{i+1}=(V_1,V_2,V_3)$ (addition modulo $m$), such that $u_i\in U_i$, $v_i\in V_i$ for all $i\in [3]$ with $|Q|\leq q$ and the vertices of $Q$ are disjoint from the vertices of $\bigcup V(K^i)$.
\end{definition}

Note that a $(q, s, n')$ connected tripartite cover contains a square cycle on at least $n'$ vertices (by ``winding around" each balanced complete tripartite graph and using each $P_i$ to get to the next tripartite graph).

%

\begin{lemma}\label{coverconnecting}
For all $0<\ep\ll \alpha\ll 1$ there exists $n_0$ such that if $G$ is a graph on $n\geq n_0$ vertices with $\delta_2(G)\geq (\frac{4}{3}-4\ep)n$ and $G$ is not $\alpha$-extremal, then the following statement holds. Given disjoint balanced complete tripartite subgraphs $K=(U_1,U_2,U_3)$ and $K'=(V_1,V_2,V_3)$ of $G$ with color classes of size at least 19, there exists a square path $P=u_1u_2u_3Qv_1v_2v_3$ where $u_i\in U_i$, $v_i\in V_i$ for all $i\in [3]$, such that $|Q|\leq 18$. 
\end{lemma}

\begin{proof}  First let $u_1\in U_1$ and $v_3\in V_3$.  If possible, we choose $u_2\in U_2\setminus L$ and $u_3\in U_3\setminus L$.  If this is not possible, then $U_i\subseteq L$ for some $i=2,3$, in which case we let $T$ be a triangle in $U_i$ and note that every vertex of $T$ is a neighbor of both $u_2$ and $u_3$.  Likewise, if possible, we choose $v_1\in V_1\setminus L$ and $v_2\in V_2\setminus L$.  If this is not possible, then  $V_j\subseteq L$ for some $j=1,2$, in which case we let $T'$ be a triangle in $V_j$ and note that every vertex of $T'$ is a neighbor of both $v_1$ and $v_2$.  Now applying Lemma \ref{generalconnecting} gives a square path $u_1u_2u_3Qv_1v_2v_3$ with $|Q|\leq 18$.  Note that since $Q$ uses at most 18 vertices, in the worst case, we would use at most $19$ vertices from a single color class.
\end{proof}

\begin{lemma}[Connected cover lemma]\label{connectedcover}
For all $0<\ep\ll \eta\ll \alpha\ll 1$ there exists $n_0$ and $c>0$ such that if $G$ is a graph on $n\geq n_0$ vertices with $\delta_2(G)\geq (\frac{4}{3}-2\ep)n$ and $G$ is not $\alpha$-extremal, then $G$ contains a $(18, c\log n, (1-2\eta)n)$ connected triangle cover. 
\end{lemma}

\begin{proof}
First apply Lemma \ref{lem:cover} to get a $(c'\log n, (1-\eta)n)$ tripartite cover $\cT= \{K^0,\ldots,K^{m-1}\}$ (where $c'$ is the constant coming from Lemma \ref{lem:cover}). Fix an orientation for each tripartite graph in $\cT$ and applying Lemma \ref{coverconnecting} to connect $K^i=(U_1,U_2,U_3)$ to $K^{i+1}=(V_1,V_2,V_3)$ by a square path $P_i =u_1u_2u_3Qv_1v_2v_3$ where $v_h\in V_h$ and $v'_h\in V'_h$ for all $h\in [3]$ and $|Q|\leq 18$.  We make all of the vertices of $P_i$ forbidden to be used for any later connection. If at some point in this process some $K\in \cT$ has more than $\frac{c'}{2}\log n$ forbidden vertices, we make all the vertices in $K$ forbidden.  Note that there are at most $\frac{n}{c'\log n}$ connections to be made, each one causing $24$ vertices to become forbidded -- with the additional rule that once half of the vertices from a color class are used, all vertices in that tripartite graph are made forbidden -- we have that the total number of forbidden vertices at any point in the process is at most $2\cdot 24 \frac{n}{c'\log n} < \eps n$. Thus, after every step the remaining graph still satisfies $\deg(u)+\deg(v) \geq (\frac{4}{3} - 4\ep)n$, hence we can continue to apply Lemma \ref{coverconnecting} until we construct $P_{m-1}$ from $K^{m-1}$ to $K^0$.  Finally we remove all vertices from the tripartite graphs that are part of some $P_i$ except the starting and ending triangles and rebalance the tripartite graphs by discarding arbitrary subset of vertices from larger color classes, noting that at most $18$ vertices could be removed from each color class; so at most $18\cdot \frac{n}{c'\log n} < \eps n$ in total.  We now have the desired $(18, c\log n, (1-2\eta)n)$ connected tripartite cover (with $c=c'/2$). 

\end{proof}

\subsection{Inserting the remaining vertices}\label{sec:inserting}

Finally we show that if we are given a connected tripartite cover, we can assign the remaining vertices to the tripartite graphs in such a way that they can be incorporated into a square cycle.

\begin{lemma}\label{lemma:inserting}
Let $0<\ep,c \ll \eta\ll \alpha\ll 1$ and $G$ be a graph on $n$ vertices containing a $(18, c\log n, (1-2\eta)n)$ connected tripartite cover $\mathcal{K}$ with square paths $\cP = \{P_1,\dots, P_{m}\}$.  If $n$ is sufficiently  large, $\delta_2(G)\geq (\frac{4}{3}-2\ep)n$, and $G$ is not $\alpha$-extremal, then $G$ contains a collection $\cT$ of disjoint (not necessarily balanced) complete tripartite graphs $\{T^1, \dots, T^m\}$, such that $|V(T^i)|\geq (1-\eta^{1/3})c\log n$ and $|V(\cT)|\geq (1-\sqrt{\eta})n$ together with a function $f$ from $V(G)-V(\cT)-V(\cP)$ to $\{T^1, \dots, T^m\}$ having the property that $|f^{-1}(T^i)|\leq \frac{1}{\eta^{1/3}}\frac{2\eta n}{m}$ and $V(T_i)\cup f^{-1}(T_i)$ contains a square path starting with the last edge of $P_{i-1}$ and ending with the first edge of $P_{i}$.
\end{lemma}

\begin{proof}

Let $U = V(G) - V(\mathcal{K}) - V(\cP)$ and note that $|U|\leq 2\eta n$.  We will try to assign the vertices of $U$ to the complete tripartite graphs, but in the process we will end up having to modify the original cover.  For convenience, we let the original cover consist of complete tripartite graphs $\{T^1,\dots, T^m\}$ and square paths $\{P_1,\dots, P_m\}$ where $\frac{1}{6c}\frac{n}{\log n}\leq m\leq  \frac{1}{3c}\frac{n}{\log n}$ and throughout the process, we will refer to the tripartite graphs by these same names even if they are modified. We assume that size of a color class in $T^i$ is $t$. However, we will maintain a set $\cT^*$ of triangles which cannot be modified as they are being used to insert vertices into some $T^i$. 


Let $w\in U$.  We will prove that we can assign $w$ to some $T^i$ while only adding at most $8$ triangles to $\cT^*$.   Once $\eta^{1/3}c\log n$ vertices have been assigned to $T^i$, then we make all of the vertices of $T^i$ forbidden.  Since there are at most $2\eta n$ vertices to be assigned, this will make at most $\frac{2\eta n}{\eta^{1/3}c\log n}$ tripartite graphs forbidden and a total of at most $\frac{2\eta n}{\eta^{1/3}c\log n}\cdot 6c\log n \leq 12\eta^{2/3} n$ forbidden vertices $Z$.  For any vertex we only consider its neighborhood in $V(G)-V(\cT^*)-V(\cP)-Z$ so for the rest of the proof we will assume that 
\begin{equation}
\delta_2(G) \geq (\frac{4}{3} - 2\ep)n- 2(|V(\cP)|-|V(\cT^*)|-|Z|)>(\frac{4}{3} - 2\sqrt{\eta})n.
\end{equation}

First, if $w$ has at least $2$ neighbors in every color class of $T^i = (T_1^i,T_2^i,T_3^i)$, then there are two triangles $(x_1,x_2,x_3)$ and $(y_1,y_2,y_3)$ in $N(w)$, such that $x_j,y_j \in T^i_j$. Clearly we can assign $w$ to $T^i$. We add the two triangles $(x_1,x_2,x_3)$ and $(y_1,y_2,y_3)$ to $\cT^*$; we say that these triangles are blocked by $w$. 

Suppose this is not the case; without loss of generality, for all $i\in [m]$ assume that $w$ has at most one neighbor in $T^i_3$.  Let $R(w)=\{v\in T^i_3\setminus N(w): |N(w)\cap T^i_1|, |N(w)\cap T^i_2|\geq 2\}$.  Since $\delta(G)\geq (\frac{1}{3}+\alpha)n>mt+\sqrt{\eta} n$, $R(w)$ is non-empty.  If $\deg(w)\leq (\frac{2}{3}-\sqrt{\eta})n$, then let $w'\in R(w)$.  By the degree condition, $\deg(w')\geq (\frac{2}{3}-\sqrt{\eta})n$.  So we may insert $w$ into $T^i$ adding two triangles to $\cT^*$ and try to insert $w'$ instead.  So we assume $\deg(w)\geq (\frac{2}{3}-\sqrt{\eta})n$ and we will try to insert $w$ by adding at most six triangles to $\cT^*$.  We may also assume that for all $v\in \{w\}\cup R(w)$, $\deg(v)<(\frac{2}{3}+\sqrt{\eta})n$, otherwise $v$ will have at least two neighbors in each color class of some $T^i$, in which case we can move $v$ to $T^i$ and replace $v$ with $w$.  This implies that for all $v\in V(G)$, if there exists $u\in \{w\}\cup R(w)$ such that $v\not\sim u$, then 
\begin{equation}\label{relaxedlower}
\deg(v)\geq (\frac{4}{3}-2\sqrt{\eta})n-(\frac{2}{3}+\sqrt{\eta})n\geq (\frac{2}{3}-3\sqrt{\eta})n.
\end{equation}
Since $\deg(w)\geq (\frac{2}{3}-\sqrt{\eta})n$ and $\deg(w, T^i_3)\leq 1$ for all $i$,  it is the case that $w$ has at least $(1-2\sqrt[4]{\eta})t$ neighbors in two color classes of at least $(1-2\sqrt[4]{\eta})m$ tripartite graphs in $\cT$, as otherwise we would have
$$  (\frac{2}{3}-\sqrt{\eta})n\leq \deg(w)\leq 2t(1-2\sqrt[4]{\eta})m+2\sqrt[4]{\eta}m(2-2\sqrt[4]{\eta})t\leq  (2-4\sqrt{\eta})tm\leq  (2-4\sqrt{\eta})\frac{n}{3}$$
a contradiction.

Let $I=\{i\in [m]: |N(w)\cap T^i_1|, |N(w)\cap T^i_2|\geq (1-2\sqrt[4]{\eta})t\}$ and $R'(w)=\bigcup_{i\in I} T^i_3\setminus N(w)$.  Note that $|R'(w)|\geq (1-2\sqrt[4]{\eta})tm\geq (\frac{1}{3}-\alpha)n$ and since $G$ is not $\alpha$-extremal, $e(R'(w))\geq \alpha n^2$.  At least $(\alpha-\sqrt{\eta})n^2$ of these edges are not incident with a triangle in $\cT^*$.  If any of these edges lie inside some $T^i_3$, then we can insert $w$ into $T^i$ (adding $w$ to $T^i_3$ unbalances $T^i_3$, but we can use the edge inside $T^i_3$ to rebalance).  Let $\alpha'=\alpha-\sqrt{\eta}$, so there are at least $\alpha' m^2$ pairs $\{i,j\}$ such that $e(T^i_3, T^j_3)\geq \alpha' t^2$ and $i,j\in I$; let $I_2$ be the set of such pairs.
\begin{claim}
Either we can insert $w$ according to the rules above or for all $\{i,j\}\in I_2$, there exists $h\in [3]$ such that (i) $e(T^i_h, T^j)<(2-\alpha'/8)t^2$, (ii) $e(T^i_h, V(G))\geq (\frac{2}{3}-3\sqrt{\eta})nt$, and (iii) $e(T^i_h, T)\leq 2t^2$ for all $T\in \cT$.  
\end{claim}

Suppose the claim holds and we are unable to insert $w$.  For each pair in $I_2$, there is some color class $T^i_h$ having the stated property.  Since there are at least $\alpha' m^2$ pairs in $I_2$ and at most $3m$ color classes, some color class $T^i_h$ has the property for at least $\alpha' m/3$ pairs.  This implies that $$\left(\frac{2}{3}-3\sqrt{\eta}\right)nt\leq e(T^i_h, V(G))\leq 2t^2m-\frac{\alpha' m}{3}\cdot\frac{\alpha' t^2}{8}\leq \left(\frac{2}{3}-\frac{\alpha'^2}{72}\right)nt,$$ a contradiction since $\eta <\alpha'$.  We now finish the proof by proving the claim.

\begin{proof}
Let $\{i,j\}\in I_2$ and let $X_i=\{x\in T^i_3: \deg(x, T^j_3)\geq \alpha' t\}$ and $X_j=\{x\in T^j_3: \deg(x, T^i_3)\geq \alpha' t\}$.  Since $e(T^i_3, T^j_3)\geq \alpha' t^2$, we have $|X_i|, |X_j|\geq \alpha' t$.  If any vertex $v\in T^i_3\cup T^j_3$ has at least $2$ neighbors in every color class of some $T\in \cT$ ($T\neq T^i$ and $T\neq T^j$), then we may move $v$ to $T$ without unbalancing and replace $v$ with $w$.  So in particular for all $v\in T^i_3$, we have
\begin{equation}\label{smallone}
\deg(v, T^j_1)\leq 1 ~\text{ or }~ \deg(v, T^j_2)\leq 1
\end{equation}

Let $X_i''=\{x\in X_i: \deg(x, T^j_1),\deg(x,T^j_2)\leq \alpha' t\}$.  If $|X_i''|\geq |X_i|/2\geq \alpha' t/2$, then $$e(T^i_3,T^j)\leq (1-\alpha'/2)t\cdot 2t+\alpha' t/2\cdot(1+2\alpha')t\leq(2-\alpha'/2)t^2$$ (clearly conditions (i) and (ii) are satisfied since $T^i_3\subseteq R(w)$); so suppose $|X_i''|<|X_i|/2$.  Since every vertex in $X_i\setminus X_i''$ has at least $\alpha' t$ neighbors in either $T^j_1$ or $T^j_2$, we can set $X_i'$ to be those vertices with at least $\alpha' t$ neighbors in $T^j_2$ and without loss of generality we may suppose $|X_i'|\geq |X_i\setminus X_i''|/2\geq \alpha' t/4$. We will now show that $T^j_1$ satisfies the conditions of the Claim.

Note that every vertex in $X_i'$ has at most one neighbor in $T^j_1$ by \eqref{smallone}, and thus 
\begin{equation}\label{incidentXi}
e(T^j_1,T^i_3)\leq t(t-|X_i'|)+|X_i'|\leq (1-\alpha'/6)t^2
\end{equation}

Also, if some vertex $w'$ in $T^j_1\cap N(w)$ has at least $2$ neighbors in each of $N(w)\cap T^i_1$ and $N(w)\cap T^i_2$, then we may move $w$ and $w'$ into $T^i_3$ and replace $w'$ with some vertex $x\in X_i'$; so suppose not.   This implies 
$$e(T^j_1, T^i_1\cup T^i_2)\leq (1-2\sqrt[4]{\eta})t\cdot t+2\sqrt[4]{\eta}t\cdot 2t=(1+2\sqrt[4]{\eta})t^2. $$  
Combining this with \eqref{incidentXi} gives
\begin{equation}
e(T^j_1, T_i)\leq (1-\alpha'/6)t^2+(1+2\sqrt[4]{\eta})t^2\leq (2-\alpha'/8)t^2
\end{equation}

If there were more than one vertex in $T^j_1$ which is adjacent to every vertex in $X_i'$, then we violate \eqref{incidentXi}; so suppose not. 
Thus by \eqref{relaxedlower}, $e(T^j_1, V(G))\geq (\frac{2}{3}-3\sqrt{\eta})nt$.  Finally if some vertex $v\in T^j_1$ has at least two neighbors in every color class of some $T\in \cT$, then we could move $v$ to $T$, replace it with a vertex from $X_i'$ (which has at least $\alpha' t$ neighbors in $T^j_2$ and $T^j_3$) and replace the vertex from $X_i'$ with $w$; thus $e(T^j_1, T)\leq 2t^2$ for all $T\neq T^i$.
\end{proof}

\end{proof}

\begin{proof}[of Theorem \ref{thm:nonextremal}] Given $G$ we first apply Lemma \ref{connectedcover} to get a $(18, c\log n, (1-2\eta)n)$ connected cover in $G$. We insert the remaining vertices into the cover using Lemma \ref{lemma:inserting} and get a set $\cT$ of $m$ disjoint complete tripartite graphs, a set of square paths $\cP$, and the function $f:V(G)-V(\cT)-V(\cP) \rightarrow \cT$. Note that for any $w$ such that $f(w)=T^i$, there are two triangles blocked by $w$. Let $(x_1,x_2,x_3)$ and $(y_1,y_2,y_3)$ be the triangles blocked by $w$, notice that by construction $x_i,y_i \in N(w)$. Create an auxiliary triangle $(z_1,z_2,z_3)$ in $T^i$ to replace these two triangles and connect $z_i$ to the common neighbors of $x_i$ and $y_i$. Note that the modified $T^i$ is still a complete tripartite graph. We similarly introduce such an auxiliary triangle for each vertex $w \in V(G)-V(\cT)-V(\cP)$. Find a triangle cover in the remaining part of $T^i$ except for the two triangles that are part of $P_{i-1}$ and $P_{i}$ by an application of Lemma \ref{cleanExtCaseLemma1}. Combining these triangles with the auxiliary triangles, we find a Hamiltonian square path by applying Lemma \ref{cleanExtCaseLemma2} that starts with the last triangle $t_{i-1}$ of $P_{i-1}$ and ends with the first triangle $t_i$ of $P_i$. 
\end{proof}

\noindent
\textbf{Acknowledgements:} We would like to thank the anonymous referees for suggesting a nicer presentation of Lemma \ref{generalconnecting}.

\end{document}